%% file: agt-5-15.tex
\documentclass{gtart_h}
\input agtout

\lognumber{15}
\volumenumber{5}
\volumeyear{2005}
\papernumber{15}
\pagenumbers{355}{368}
\received{22 January 2005} 
\revised{30 March 2005}
\accepted{12 April 2005}
\published{21 April 2005}

\usepackage{amsmath}
\usepackage[all]{xy} %This is XY-pic.
\usepackage{amssymb}

\newcommand{\ra}{\rightarrow}

\newcommand{\BR}{\mathbb{R}}

\newcommand{\BZ}{\mathbb{Z}}

\newcommand{\Kahler}{K\"{a}hler}

\newtheorem{theorem}{Theorem}
\newtheorem{thm}[theorem]{Theorem}
\newtheorem{lemma}[theorem]{Lemma}
\newtheorem{lem}[theorem]{Lemma}

\newtheorem{proposition}[theorem]{Proposition}
\theoremstyle{definition}

\newtheorem{defn}[theorem]{Definition}

\newtheorem{rem}[theorem]{Remark}
\newtheorem{quest}[theorem]{Question}

\begin{document}

\title{Geography of symplectic 4--manifolds\\with Kodaira dimension one}
\asciititle{Geography of symplectic 4-manifolds with Kodaira dimension one}

\authors{Scott Baldridge\\Tian-Jun Li}

\address{Department of Mathematics,  Louisiana State University\\Baton Rouge, 
LA 70803, USA}
\secondaddress{School of Mathematics, University of 
Minnesota\\Minneapolis, MN 55455, USA}

\asciiaddress{Department of Mathematics,  Louisiana State University\\Baton Rouge, 
LA 70803, USA\\and\\School of Mathematics, University of 
Minnesota\\Minneapolis, MN 55455, USA}

\gtemail{\mailto{sbaldrid@math.lsu.edu}{\rm\qua and\qua}\mailto{tjli@math.umn.edu}}
\asciiemail{sbaldrid@math.lsu.edu, tjli@math.umn.edu}

\begin{abstract}
The geography problem is usually stated for simply connected
symplectic 4--manifolds.  When the first cohomology is nontrivial,
however, one can restate the problem taking into account how close the
symplectic manifold is to satisfying the conclusion of the Hard
Lefschetz Theorem, which is measured by a nonnegative integer called
the degeneracy.  In this paper we include the degeneracy as an extra
parameter in the geography problem and show how to fill out the
geography of symplectic 4--manifolds with Kodaira dimension 1 for all
admissible triples.
\end{abstract}

\asciiabstract%
{The geography problem is usually stated for simply connected
symplectic 4-manifolds.  When the first cohomology is nontrivial,
however, one can restate the problem taking into account how close the
symplectic manifold is to satisfying the conclusion of the Hard
Lefschetz Theorem, which is measured by a nonnegative integer called
the degeneracy.  In this paper we include the degeneracy as an extra
parameter in the geography problem and show how to fill out the
geography of symplectic 4-manifolds with Kodaira dimension 1 for all
admissible triples.}


\maketitle

\section{Introduction}

For a minimal symplectic $4$--manifold $M$ with symplectic form
$\omega$ and symplectic canonical class $K_{\omega}$, the Kodaira
dimension of $(M,\omega)$
 is defined in the following way:
$$\kappa(M,\omega) = \left\{ \begin{array}{cl} \ \  -\infty \ \  &
\mbox{if\qua $K_\omega^2<0$\qua or \qua\ $K_\omega\cdot
[\omega]<0$}\\
0 & \mbox{if\qua $K_{\omega}^2=0$\qua and\qua $K_\omega\cdot [\omega] =0$}\\
1 & \mbox{if\qua $K_{\omega}^2=0$\qua and\qua $K_\omega\cdot [\omega] >0$}\\
2 & \mbox{if\qua $K_{\omega}^2>0$\qua and\qua $K_\omega\cdot [\omega]
>0$.}\end{array}\right.
$$
The Kodaira dimension of a non-minimal manifold is defined to be
that of any of its minimal models (see 
\cite{symp:mcduff_salamon:intro_to_symp_top}, \cite{symp:li:kod_dim_zero}).

Minimal symplectic manifolds of Kodaira dimension $-\infty$ were
classified in \cite{symp:app_gen_wall_cross}. Such manifolds are  either ${\bf CP}^2$ or an $S^2$--bundle
over a surface. Minimal symplectic manifolds of Kodaira dimension
zero were studied in \cite{symp:li:kod_dim_zero}: it was speculated
that they are  either K3, Enriques surface or a $T^2$--bundle over
$T^2$; and it was shown that the $\chi$ and $\sigma$ are bounded
if $b_1$ is bounded by $4$.

For symplectic $4$--manifolds with Kodaira dimension 1 or 2 we cannot
expect to have a classification: Gompf \cite{symp:gompf:new_const}
showed that any finitely presented group is the fundamental group
of a symplectic $4$--manifold either of Kodaira dimension 1 or of
Kodaira dimension 2. Instead one is interested in further
illustrating the diversity of simply connected minimal symplectic
$4$--manifolds (\cite{symp:mccarthy:normal_connect_sum}, \cite{symp:gompf:new_const},
\cite{sw:knots_links_and_four_man}, \cite{symp:stip:example}). The
problem is to realize all pairs  $(\chi(M), \sigma(M))$  of a
simply connected minimal symplectic $4$--manifold $M$ subject to the
Noether condition $2\chi+3\sigma \equiv \sigma\pmod 8$ (symplectic
manifolds admit almost complex structures) and the conjectured
inequality $\chi\geq 3\sigma$ (the symplectic
Bogomolov-Miyaoka-Yau inequality).

For simply connected minimal symplectic $4$--manifolds with Kodaira
dimension 1 this has a simple positive answer. Such a manifold has
$b_1=0$ and $K^2=0$, so
\begin{equation}
2\chi+3\sigma=4(1-b_1+b_+)+ \sigma=0.
\label{Eq:Kod=1}\end{equation}
Therefore $\sigma$ is nonpositive and $\chi$ is determined by
$\sigma$. And since $M$ is almost complex, $b^+-b_1$ is odd, and hence
$\sigma(M)$ is divisible by 8.

The Dolgachev surfaces and Elliptic surfaces $E(n)$ with $n\geq 2$ are
simply connected minimal symplectic $4$--manifolds with Kodaira
dimension 1 with signature $-8$ and  $-8n$ respectively.
Hence we have the well-known fact:

\begin{proposition} Every negative integer divisible by $8$ is the signature of a
simply connected minimal symplectic $4$--manifold with Kodaira
dimension 1. \label{prop:simply_connected_examples}
\end{proposition}

In this paper we  investigate the geography question of minimal
symplectic $4$--manifolds of Kodaira dimension 1 by taking into account
the first Betti number $b_1$ and the  cup product structure on $H^1$
(a question perhaps mostly of interest to $4$--manifold topologists).
By Equation \eqref{Eq:Kod=1}, the pair $(\chi,\sigma)$ is equivalent
to $(\sigma, b_1)$ when the Kodaira dimension is 1. We will use the
latter pair of numbers in what follows to make definitions and
statements  more transparent.  The third parameter is determined by
the symplectic form and can be formulated in terms of a \Kahler-like
condition called Lefschetz type.

\begin{defn}
Symplectic $4$--manifolds $(M,\omega)$ are said to be of {\em
Lefschetz type} if  ${[\omega]\in H^2(M;\BR)}$ satisfies the
conclusion of the Hard Lefschetz Theorem, namely, that $\cup
[\omega]:H^1(M;\BR) \ra H^3(M;\BR)$  is an isomorphism.
\end{defn}

Based upon that definition one gets:

\begin{defn}
The {\em degeneracy} $d(M,\omega)$ is the rank of the kernel of
the map
$$\cup [\omega]: H^1(M;\BR) \ra H^3(M;\BR).$$
\end{defn}

The first example of a symplectic $4$--manifold with  nonzero
degeneracy occurs when $b_1=2$ and $b_+=1$. A symplectic $4$---manifold
with $b_+=1$ and $\kappa=1$ must satisfy $b_1=0$ or $b_1=2$ using
the fact that $2\chi +3\sigma=0$. Examples of $b_1{=}2$ 4--manifolds
of Lefschetz type are easy to construct and have been known for some
time, but examples of manifolds not of Lefschetz type remained
unknown until Baldridge \cite{symp:baldridge:newsymp_4man}. These
$b_+=1$ manifolds with nonzero degeneracy are the starting point for
the examples described in this paper.

\begin{defn} Any triple $(a,b,c) \in \BZ^3$ is called an {\em admissible
triple} {\em(or Lefschetz admissible)} if  $a=8k$ where $k$ is a
non-positive integer, $0\leq c \leq b$, $b-c$ even, and $b \geq
\max\{0, 2 + a/4\}$.
\end{defn}

In the above definition, $a$ corresponds to the signature of a
symplectic 4--manifold, $b$ is the first Betti number, and $c$ is
the degeneracy.  This definition covers all triples except possible
counterexamples to the conjectured BMY inequality. The next lemma
explains why $b-c$ should be even.

\begin{lemma}
Let $(M,\omega)$ be a closed symplectic 4--manifold.  The
skew-symmetric bilinear form $Q_M:H^1(M;\BR)\times H^1(M;\BR) \ra
\BR$ defined by
$$Q_M(a,b) = \int_M a\cup b\cup [\omega]$$
has rank $b_1(M)-d(M,\omega)$.  Furthermore, $rank \ Q_M$ is even.
\end{lemma}

\begin{proof}
Pick a compatible metric $g$ and set $$K= \{\alpha\in H^1(M;\BR)
\; | \; \alpha\cup [\omega] =0 \in H^3(M;\BR)\}.$$  $K$ is a closed
subspace of $H^1(M;\BR)$ with rank $d(M,\omega)$. Let $W$ be the
orthogonal complement of $K$ in $H^1(M;\BR)$ with respect to the
$L^2$--norm on $H^1(M;\BR)$.

We claim that $Q_M$ is nondegenerate on $W$.  To see this, suppose
$a \in W$ satisfies $Q_M(a,b)=0$ for all $b\in W$.  If the class
$a\cup [\omega]$ is not zero, then by Poincar\'{e} duality there
exists a $1$-cocycle $\gamma$ such that $Q_M(\gamma,a) = \langle
\gamma\cup a \cup [\omega], [M]\rangle \not=0$. Therefore
$a\cup[\omega] = 0$ and $a\in K$ which implies $a=0$.

Now by the same proof that symplectic forms are locally isomorphic
to the standard symplectic form on $(\BR^{2n}, \omega_0)$, $rank
\; Q_M$ is even dimensional.
\end{proof}

We can now state the main theorem of this paper:

\begin{thm} \label{thm:main_thm} For any admissible triple $(a,b,c)$  there exists a minimal
symplectic 4--manifold  $(M,\omega)$ of Kodaira dimension
$\kappa(M)=1$ with
$$(a,b,c)=(\sigma(M),b_1(M), d(M,\omega)),$$
where $\sigma(M)$ is the signature and $b_1(M)$ is the first Betti
number.
\end{thm}

\section{An important lemma}
\label{sec:symp_const}

In Section \ref{sec:examples} we will produce examples of
symplectic 4--manifolds for admissible triples $(0,b,c)$. We will
then use these manifolds to fill out triples where $a<0$ in
Section \ref{sec:prove_main_theorem}. To build these manifolds we
need a general theorem that was first reported in
\cite{symp:sympmanfreecircleaction}.

\begin{thm} \label{thm:symp_constr}Let $\Sigma$ be a closed,
oriented, connected surface of genus $g$ with an
orientation preserving diffeomorphism $\varphi:\Sigma \ra \Sigma$.
Let $Y$ be the mapping torus with respect to $\Sigma$ and
$\varphi$ and $p:Y\ra S^1$.  Let $\mathfrak{e} \in H^2(Y;\BZ)$ be
any class such that $\mathfrak{e}|_{\Sigma} = 0$ and let $M \ra Y$
be the $S^1$-bundle over $Y$ with Euler class $\mathfrak{e}$. Then
$M$ is a smooth oriented closed symplectic 4--manifold.
\end{thm}

The advantage of the construction below over the one in
\cite{symp:sympmanfreecircleaction} is that we can derive the
cohomology ring explicitly for the examples we need.

\proof
We begin by constructing a basis for the cohomology $H^1(Y;\BZ)$
from smooth integral 1--forms.  Let $\theta \in \Omega^1(Y)$ be
the pullback of the volume of $S^1$. This is a closed, integral,
nonzero 1--form.

Consider the map $\varphi^* -1$ on $H^1(\Sigma;\BZ)$ and let $k=
\mbox{rank} \ker (\varphi^*-1)$. The Wang exact sequence implies
that the rank of $H^1(Y)$ is $k+1$:
$$
\xymatrix{ H^0(\Sigma) \ar[r] \ar@{=}[d] & H^1(Y) \ar@{=}[d]
\ar[r]^{\ |_{\Sigma}} &H^1(\Sigma) \ar@{=}[d] \ar[r]^{\varphi^*-1}
& H^1(\Sigma)
\ar@{=}[d] \ar[r]^{\mu} & H^2(Y) \ar[d] \\
\BZ \ar[r] & \BZ\oplus\BZ^k  \ar[r] & \BZ^k \oplus
\BZ^{2g-2-k}\ar[r] & \BZ^{2g-2-k} \oplus \BZ^k \ar[r] & \BZ^k
\oplus \BZ. }
$$
Write down a basis for $H^1(Y;\BZ)$ as follows.  Let
$\langle \gamma_1,\cdots, \gamma_k\rangle$ be a basis of closed
integral 1--forms on $\Sigma$ for the subspace of
$H^1(\Sigma;\BZ)$ which is preserved by $\varphi^*$. Extend this
to a basis of $H^1(\Sigma;\BZ)$ by closed 1--forms:
$$\langle [\gamma_1],\dots, [\gamma_k], [\epsilon_{k+1}], \ldots,
[\epsilon_{2g}]\rangle.$$
Because $[\gamma_i]$ is invariant under $\varphi^*$ there is a
function $f_i\in \Omega^0(\Sigma)$ such that $\varphi^*(\gamma_i)
= \gamma_i + df_i$ point-wise. To construct a closed 1--form on
$Y$ which is non-trivial in cohomology, choose a smooth function
$\rho:[0,1] \ra [0,1]$ which is identically 0 near 0 and
identically 1 near 1 and extend $\gamma_i$ to $\Sigma \times
[0,1]$ by defining
$$\bar{\gamma}_i(x,t) = \rho(t)\gamma_i(x) + 
\big(1-\rho(t)\big)\big(\varphi^*(\gamma_i(x))\big)
- \Big(\frac{d}{dt}\rho(t)\Big)f_i(x)dt.$$ Since $\bar{\gamma}_i$ can be
identified on the boundary of $\Sigma \times [0,1]$ using
$\varphi^*$, and
$$\frac{d^n}{dt^n}\varphi^*(\bar{\gamma}_i)|_{\Sigma\times\{0\}} =
\frac{d^n}{dt^n}\bar{\gamma}_i|_{\Sigma\times\{1\}}$$ for all
whole numbers $n$, $\bar{\gamma}_i$ is a closed smooth 1--form on
$Y$. Hence,
$$H^1(Y;\BR) = \langle [\theta], [\bar{\gamma}_1], \ldots,
[\bar{\gamma}_k]\rangle.$$
Next we construct a basis of $H^2(Y)$. First, there is a smooth
closed integral 2--form $\Omega_\Sigma \in \Omega^2(Y)$ which
restricts to the volume form of $\Sigma$ on each fiber of $p:Y\ra
S^1$. Writing down a basis which spans $\mbox{Im}
(\mu:H^1(\Sigma)\ra H^2(Y))$ in the Wang sequence is crucial to
understanding what symplectic 4--manifolds can be built from $Y$.
Note that
$$\varphi_*(PD(\gamma_i)) = \varphi_*(\gamma_i\cap [\Sigma]) =
\varphi_*(\varphi^*(\gamma_i))\cap [\Sigma]) = \gamma_i \cap
\varphi_*[\Sigma] = PD(\gamma_i),$$ where $PD$ is the Poincar\'{e}
dual.  Hence  $PD(\gamma_i)$ is $\varphi$-invariant in homology
and $\langle t, PD(\gamma_1), \ldots, PD(\gamma_k)\rangle$ is a
basis for $H_1(Y)$ where $t$ is the loop $pt \times [0,1]$ in
$\Sigma\times [0,1]$ for some fixed point of $\varphi$.

Let $\xi_i$ be a 1--form which is the Hom-dual of $PD(\gamma_i)$.
Extend the linearly independent set $\langle
[\xi_1],\dots,[\xi_k]\rangle$ to a basis of $H^1(\Sigma)$ given by
$$\langle [\xi_1],\dots, [\xi_k], [\zeta_{k+1}], \ldots
[\zeta_{2g}]\rangle.$$ The classes $[\xi_i]$ are not necessarily
$\varphi^*$-invariant. Nevertheless, since
$$\langle \varphi^*(\xi_i), PD(\gamma_j) \rangle = \langle \xi_i,
\varphi_*(PD(\gamma_j))\rangle = \langle \xi_i,PD(\gamma_j)\rangle
=\delta_{ij},$$ there exists functions $g_i \in \Omega^0(\Sigma)$
and integers $c_j$ such that
\begin{eqnarray}\varphi^*(\xi_i)  = \xi_i +
\sum_{j=k+1}^{2g} c_j \zeta_j + dg_i\label{eqn:varphi_on_xi}
\end{eqnarray}
point-wise. Using $\xi_i$ and $g_i$, construct a smooth 1--form on
$Y$ by specifying
$$\bar{\xi}_i(x,t) = \rho(t)\xi_i(x) + \big(1-\rho(t)\big)\big(\varphi^*(\xi_i(x))\big)
- \Big(\frac{d}{dt}\rho(t)\Big)g_i(x)dt.$$ This 1--form is not necessarily
closed, in fact,
$$d\bar{\xi}_i =  \Big(\sum c_j\zeta_j\Big) \wedge
\Big(\frac{d}{dt}{\rho}\Big) \theta.$$ However $\bar{\xi}_i\wedge \theta$ is a
smooth closed integral 2--form and $$\langle [\Omega_\Sigma],
[\bar{\xi}_1\wedge \theta], \ldots,
[\bar{\xi}_k\wedge\theta]\rangle$$ is a basis for $H^2(Y;\BZ)$.

We can now complete the proof of Theorem \ref{thm:symp_constr}. If
$\mathfrak{e} \in H^2(Y;\BZ)$ such that
$\mathfrak{e}|_{\Sigma}=0$, then $$\mathfrak{e} = \sum_{i=1}^k p_i
[\bar{\xi}_i\wedge \theta]$$ where $p_i\in \BZ$. Let $M$ be the
$S^1$-bundle over $Y$ with Euler class $\mathfrak{e}$ and
connection 1-form $\eta$ given by $d\eta =\sum_{i=1}^k p_i
\pi^*(\bar{\xi}_i\wedge \theta)$ point-wise. Then
$$\omega = \pi^*(\Omega) + \pi^*(\theta)\wedge \eta$$ is an everywhere
nondegenerate 2--form on $M$ such that
$$d\omega =-\pi^*(\theta)\wedge d\eta = -\pi^*(\theta) \wedge \sum_{i=1}^k p_i \pi^*(\bar{\xi}_i\wedge
\theta) = 0.$$ This ends the proof of the theorem.\endproof

McDuff and Salamon \cite{symp:mcduff_salamon:intro_to_symp_top}
raised the question whether there is a free symplectic circle
action on symplectic $4$--manifold with contractible orbits. They
pointed out that null-homologous orbits exists on certain
$T^2$--bundles over $T^2$. Here we show that there cannot be any
contractible orbits.

\begin{proposition} Suppose $S^1$ acts freely (fixed point free)
 on $(M,\omega)$ preserving
$\omega$. Then the orbits are essential in $\pi_1(M)$.%
\footnote{Added in proof: We have learned that this is a special case
of a result of Kotschick \cite[Theorem
1]{symp:kotschick:free_circ_act}, where the action is not assumed to
preserve the symplectic structure in dimension 4.  Moreover examples
of symplectic free actions with contractible orbits in every dimension
at least 6 are constructed in \cite[Theorem
2]{symp:kotschick:free_circ_act}.}
\end{proposition}

\begin{proof} Let $Y=M/S^1$. Then it is known that $Y$ fibers over $S^1$.
When the fiber is of genus at least one, $Y$ is a $K(\pi, 1)$
space from the homotopy exact sequence associated to this
fibration, in particular it has trivial $\pi_2$. Consider the
homotopy exact sequence associated to the fibration
$S^1\longrightarrow M\longrightarrow Y$,
$$\cdots\longrightarrow
\pi_2(Y)\longrightarrow \pi_1(S^1)\longrightarrow \pi_1(M)
\longrightarrow\cdots.$$ The orbits are contractible implies that
$\pi_2(Y)$ surjects onto $\pi_1(S^1)={\bf Z}$. This is a
contradiction.
\end{proof}

\section{Bundle manifolds}
\label{sec:bundleman}

In this section we describe symplectic 4--manifolds which  are
special $S^1$-bundles over a base which is itself a surface bundle
over $S^1$, which we will call {\em bundle manifolds}. They are
uniquely specified by four whole numbers: three weights $g,d,k$
where $0\leq d \leq k \leq g$ and a number $e=0,1,2$ based upon
the Euler class of the $S^1$-bundle. We denote these manifolds by
$B(d,k,g;e)$.

Construct a bundle manifold as follows. Let $\Sigma$ be a surface
of genus $g>0$. Let $\langle a_1,b_1,\ldots a_g,b_g\rangle$ be
smooth loops which represent the usual symplectic basis of $H_1(\Sigma)$ such that $a_i\cdot b_j= \delta_{ij}$ and zero otherwise.

%\begin{figure}[h]
%\begin{center}
%\psfrag{a1}{$a_1$} \psfrag{b1}{$b_1$}\psfrag{a2}{$a_2$}
%   \includegraphics{surface.eps}
%    \caption{Symplectic basis on a genus 3 surface. }
%    \label{fig:surface_genus_g}
%\end{center}
%\end{figure}

Similarly, let $\langle \alpha_1, \beta_1, \ldots \alpha_g,
\beta_g\rangle$ be closed integral 1--forms which represent a dual
basis with respect to the $a_i$'s and $b_i$'s:
$$\alpha_i(a_j) = \delta_{ij}, \hspace{1cm} \beta_i(b_j) = \delta_{ij}, \hspace{1cm} \mbox{ and
zero otherwise.}$$
Consider the diffeomorphism given by the following sequence of
Dehn twists acting on the left,
$$
\varphi = (T_{b_g}T^{-1}_{a_g}) \cdots
(T_{b_{k+1}}T^{-1}_{a_{k+1}}) \cdot T_{a_d}\cdots T_{a_1}.
$$
This diffeomorphism has the following properties:

\begin{itemize}
\item  For each $1\leq i \leq d$ when $d\not=0$,  the subspace spanned by
$\langle \alpha_i,\beta_i\rangle$ has a 1--dimensional subspace
preserved by $\varphi^*$,
$$\varphi^*\alpha_i = \alpha_i+ \beta_i \mbox{ \ \ and \ \ }
\varphi^*\beta_i = \beta_i$$ as cohomology classes.

\item For each $d< i \leq k$ when $d\not=k$, the subspace spanned by
$\langle\alpha_i,\beta_i\rangle$ is preserved by $\varphi^*:H^1(\Sigma;\BZ)\ra H^1(\Sigma;\BZ)$.

\item The subspace spanned by $\langle  \alpha_{k+1}, \beta_{k+1},
\ldots \alpha_{g}, \beta_{g}\rangle$ contains no subspace which is
preserved by $\varphi^*$.
\end{itemize}

The mapping torus
$$Y = \left(\Sigma \times [0,1]\right) / \left((x,1) \sim (\varphi(x),0)\right)$$
is a smooth, closed, 3--dimensional manifold which, by the proof
of Theorem~\ref{thm:symp_constr}, has the following basis in
cohomology:
\begin{eqnarray} \label{eq:explicit_basis}
H^1(Y) & = & \langle \theta,  \ \ \bar{\beta}_1, \ldots,
\bar{\beta_d}, \ \ \bar{\alpha}_{d+1}, \bar{\beta}_{d+1}, \ldots,
\bar{\alpha}_k,\bar{\beta}_{k}\rangle \\ H^2(Y) & = & \langle
\Omega, \ \ \bar{\alpha}_1\wedge \theta, \ldots,
\bar{\alpha}_d\wedge \theta, \nonumber \\ && \hspace{1cm}
\bar{\alpha}_{d+1}\wedge \theta,\bar{\beta}_{d+1}\wedge\theta,
\ldots,
\bar{\alpha}_k\wedge\theta,\bar{\beta}_k\wedge\theta\rangle.
\nonumber
\end{eqnarray}
Also by Theorem \ref{thm:symp_constr}, any $S^1$-bundle over $Y$
is symplectic as long as the Euler class $\mathfrak{e} \in
H^2(Y;\BZ)$ is zero when restricted to a fiber $\Sigma$. Therefore
we can specify a symplectic 4--manifold $B(d,k,g;e)$ using the
diffeomorphism $\varphi$ by choosing four whole numbers $g,d,k,e$
such that $0\leq d\leq k\leq g$ and
$$e =
\left\{\begin{array}{c@{\hspace{.2cm}}l}
0 & \mathfrak{e} = 0,\\
1 & \mathfrak{e} = [\bar{\alpha_i}\wedge \theta] \mbox{ for some } 0 < i \leq d, d\not=0\\
2 & \mathfrak{e} \in \mbox{span}\langle [\bar{\alpha}_{d+1}\wedge
\theta],[\bar{\beta}_{d+1}\wedge\theta], \ldots,
[\bar{\beta}_k\wedge\theta]\rangle, \mbox{ } \mathfrak{e} \mbox{
primitive, and }d\not=k. \end{array}\right.$$

\begin{rem} The manifolds described in \cite{symp:baldridge:newsymp_4man} are
just the bundle manifolds given by $B(1,1,g;1)$ for all genera
$g>1$.
\end{rem}

The bundle manifold $B(d,k,g;e)$ does not depend on the
choice of an Euler class $\mathfrak{e}\in H^2(Y)$; two bundle
manifolds with different Euler classes but with the same data
$g,d,k,e$ are diffeomorphic.

The following lemmas will be helpful in the next section.

\begin{lem} \label{lem:Kodaira_dim} The Kodaira dimension of a bundle manifold
$B(d,k,g;e)$ is $\kappa=0$ if $g=1$ and $\kappa=1$ if $g>1$.
\end{lem}

\begin{proof} This theorem follows using similar arguments as in
\cite{symp:baldridge:newsymp_4man}.
\end{proof}

\begin{lem} The signature of $B(d,k,g;e)$ is
zero.
\end{lem}

\begin{proof} This is an easy computation given the existence of a free circle
action on bundle manifolds.
\end{proof}

\begin{lem} \label{lem:first_betti_num}
The first Betti number of the bundle manifold  $B(d,k,g;e)$ is:
$$b_1 = \left\{\begin{array}{c@{ \ \ \ }c}  2k-d+2 &
e=0\\2k-d+1  & e\not=0.\end{array}\right.$$
\end{lem}

\begin{proof} The Gysin Sequence,
\begin{eqnarray*}
\xymatrix{ 0 \ar[r] & H^1(Y) \ar[r]^{\pi^* \ \ \ \ \ } & H^1(B(d,k,g;e))
\ar[r] & H^0(Y) \ar[r]^{\cup \mathfrak{e}} & H^2(Y), }
\end{eqnarray*}
implies that
\begin{eqnarray*}
H^1(B(d,k,g;e),\BZ) & \cong & \left\{ \begin{array}{l@{\ \ \ }l}
 H^1(Y,\BZ) \oplus \BZ, &
e=0\\H^1(Y;\BZ), & e\not=0.\end{array}
\right. \\
\end{eqnarray*}
The first Betti number can then be calculated using the
basis constructed in Equation \eqref{eq:explicit_basis}.\end{proof} 

\begin{lem}\label{lem:degeneracy}The degeneracy $d((B(d,k,g;e),\omega)$  is
$$d(B(d,k,g;e),\omega)) = \left\{\begin{array}{c@{ \ \  \ }c}d &
e=0\\ d+1 & e\not=0.\end{array}\right.$$
\end{lem}

\begin{proof} We prove the case when $e=0$, i.e., $B(d,k,g;0) = Y \times S^1$, noting that the case when $e\not=0$ is similar.  Let $t$ be a section of the fibration $Y \ra S^1$ such that $\langle \theta,[t] \rangle=1$.  A basis for $H_3(Y\times S^1;\BZ)$ can be described as follows.  By construction, $\varphi_*(a_i)=a_i$ for $1\leq i\leq d$. While $a_i$ is null-homologous due to the relation created by $\varphi_*(b_i)=a_i+b_i$, the 2-cycle $a_i \times t$ is not, and the space generated by $\{a_i\times t\times S^1\, | \, 1\leq i\leq d\}$ forms a linearly independent subspace of  $H_3(Y\times S^1;\BZ)$.  For $d<i\leq k$, $\varphi_*$ is identity on $a_i$'s and $b_i$'s,
and so the space generated by $\{a_i\times t \times S^1, b_i\times t \times S^1 \, | \, d<i\leq k\}$ is also a linearly independent subspace.  Altogether, a basis for homology 3-cycles of $Y\times S^1$ is
\begin{eqnarray*}
H_3(Y\times S^1;\BZ) & = & \langle [Y], \Sigma\times S^1, \, a_1 \!\times\! t \!\times\! S^1,\, \ldots,\, a_d\!\times \! t \!\times \! S^1, \\ & & \hspace{.25cm} a_{d+1}\!\times\! t \! \times \! S^1, \, b_{d+1} \!\times \! t \! \times \! S^1, \, \ldots, \, b_k \! \times \! t \! \times \! S^1 \rangle.
\end{eqnarray*}
Using the basis described in Equation \eqref{eq:explicit_basis} we can determine the kernel of $\cup \omega:H^1(Y\times S^1;\BR) \ra H^3(Y\times S^1;\BR)$.  Clearly $\theta \wedge \omega = \pi^*(\theta \wedge \Omega)$ evaluates nonzero  on $[Y]$, so $\theta$ is not in the kernel.   Similarly, $\eta$ (the connection 1-form for the $S^1$ factor of $Y\times S^1$) is not in the kernel because $\eta \wedge \omega = \eta \wedge \pi^*\Omega$  evaluates nonzero on $\Sigma \times S^1$.  Observe that  $\langle \overline{\alpha}_i \wedge \omega, a_i\!\times\! t\!\times \! S^1\rangle \not= 0$ and $\langle \overline{\beta}_i \wedge \omega, b_i\!\times\! t \!\times \! S^1\rangle \not= 0$ for $d<i\leq k$.  Finally,  $\overline{\beta}_i\wedge \omega$'s for $1\leq i \leq d$ evaluates zero on all 3--cycles except possibly $a_i\!\times\! t\!\times \! S^1$.  However $\langle \overline{\beta}_i \wedge \omega, a_i\!\times\! t \!\times \! S^1\rangle = \pm \langle \beta_i, a_i\rangle = 0$.  Therefore the subspace generated by $\overline{\beta}_i$'s for $1\leq i\leq d$ is the kernel of $\cup \omega$ implying that $d(B(d,k,g;0)=d$.
\end{proof}

\section{Examples of admissible triples with $a=0$}
\label{sec:examples}

We first prove Theorem \ref{thm:main_thm} for admissible triples
$(a,b,c)$ when $a=0$ (i.e. where the signature of the manifold is
zero). In the next section we use these manifolds to construct new
manifolds for admissible triples when $a<0$.

\medskip

\noindent{\bf The case where $a=0$ and  $b$ is even}\qua\ Fix
$b=2l$ for some whole number $l>0$. We are looking for symplectic
4--manifolds $(M,\omega)$ with Kodaira dimension 1 and
$\sigma(M)=0$, $b_1(M)=2l$, and even $d(M,\omega)$ where $0 \leq
d(M,\omega) \leq 2l$. When we restrict to bundle manifolds with
$e=0$, we get by Lemma \ref{lem:first_betti_num} that
$$d=2(k-l+1).$$ Since $0\leq d\leq k$, the possible bundle manifolds occur when $k=l-1, \cdots, 2l-2$.  This,
together with Lemma \ref{lem:degeneracy}, implies that the
infinite collection of symplectic 4--manifolds
$$ \{ B_0(i):=B(2i, \ l-1 +i, \ g; \ 0) \ \ | \ \ i =0,1,\ldots, l-1 \mbox{ and } g \geq \mbox{max}(l-1+i, 2)\}$$ satisfy
\begin{eqnarray*}
b_1(B_0(i)) &\ = \ & 2l = b \mbox{\qua  and  }\\
d(B_0(i), \omega) & \  =  \ & 2i.
\end{eqnarray*}
This leaves only the case of an admissible triple $(0,b,c)$ where
$c = b$. For that we investigate the case when $e=1$. Then by
Lemma \ref{lem:first_betti_num}, bundle manifolds with $e=1$ exist, and
satisfy
$$d= 2(k-l) +1$$
when $1 \leq d\leq k$.  Hence possible values of $k=l,l+1,\ldots,
2l-1$ and the infinite set of bundle manifolds
$$\{B_1(i) := B(2i+1, \ l+i, \ g; \ 1) \ \ | \ \ i = 0,1,\ldots,
l-1 \mbox{\qua and\qua }  g\geq \mbox{max}(l+i, \ 2)\}$$ satisfy
\begin{eqnarray*}
b_1(B_1(i)) &\ = \ & 2(l+i)-(2i+1) + 1 \ = \ b \mbox{\qua  and  }\\
d(B_1(i), \omega) & \  =  \ & (2i+1) + 1 = 2(i+1)
\end{eqnarray*}
by  Lemma \ref{lem:first_betti_num} and Lemma
\ref{lem:degeneracy}. Therefore the manifold $B_1(l-1)$ is an
example of an admissible triple $(0,b,b)$. This ends the case when
$b$ is even.

\medskip
\noindent{\bf The case when $a=0$ and $b$ is odd}\qua\  Fix
$b=2l+1$ where $l> 0$.  We restrict to bundle manifolds where
$e\not=0$.  Lemma \ref{lem:first_betti_num} implies that $d$ is a
function of $k$,
$$d=2(k-l),$$
where $k$ can equal $l$ only if $e=2$, $k=2l$ only if $e=1$, and $k=l+1,
\ldots, 2l-1$ for $e=1 \mbox{ or } 2$. The infinite set of bundle
manifolds
\begin{eqnarray*}
 B_1(i) &:=& B(2i,\ l+i,\ g; \ 1)\ \ \ \  \ \ i=1,\ldots, l
\mbox{\qua and\qua } g\geq \mbox{max}(l+i,2),\\
 B_2(i) &:=& B(2i,\ l+i,\ g; \ 2)\ \ \ \  \ \ i=0,1,\ldots, l-1
\mbox{\qua and\qua } g\geq \mbox{max}(l+i,2),
\end{eqnarray*}
satisfies
\begin{eqnarray*}
b_1(B_e(i)) &\ = \ & 2(l+i)-2i+1 = b \mbox{\qua  and}\\
d(B_e(i), \omega) & \  =  \ & 2i+1
\end{eqnarray*}
for $e=1$ or $2$ by Lemma \ref{lem:first_betti_num} and Lemma
\ref{lem:degeneracy}. This completes the proof of Theorem
\ref{thm:main_thm} for the case when $a=0$.

\begin{rem} We could have generated examples for $b$ odd by looking at
bundle manifolds with $e=0$.  In that case there are examples
for $b-c$ even and $b\not=1$ except for the admissible triple
$(0,b,c)$ where $b=c$.
\end{rem}

\section{Proof of Theorem \ref{thm:main_thm}}
\label{sec:prove_main_theorem}

To fill out all admissible triples $(a,b,c)$ for $a<0$ and
divisible by $8$, we fibersum $E(n)$ with bundle manifolds. This
turns out to be much easier than the $a=0$ case since we need only
work with $B(d,k,g;e)$ where $e=0$.

Inside $B(d,k,g;0)$ we can find a symplectic torus $T$ with $T
\cdot T=0$, constructed as follows.  Since $B(d,k,g;0) = Y\times
S^1$, recall the definition of $Y$ and let $t$ be a section of
$p:Y\ra S^1$. Then $T = \pi^{-1}(t)$ is a nontrivial torus in
$B(d,k,g;0)$ which clearly has $T \cdot T  =0$. Letting $s$ be the
loop generated by the $S^1$ factor of $Y\times S^1$, we can write
$T = t\times s \in H_2(B(d,k,g;0);\BZ)$. Since $\omega|T \not=0$
point-wise, $T$ is a symplectic submanifold.

Next let $T'$ be a generic torus fiber of $E(n)$ in the neighborhood
of a cusp fiber. Consider the fibersum $$E(n, d,k,g) = E(n)\#_{T'=T}
B(d,k,g;0)$$ where $n\geq 2$.  Since we are fibering along two
symplectic surfaces, $E(n,d,k,g)$ is symplectic with symplectic form
$\tilde{\omega}$. It was proved in
\cite{symp:li-stip:minimal_connect_sum} that fiber-sums of minimal
manifolds are minimal. Therefore  all the manifolds above are
minimal.

Clearly the two loops $t$ and $s$ become null-homologous in
$H_1(E(n,d,k,g);\BZ)$ after identifying along the boundary, making
$b_1(E(n,d,k,g) = 2k-d$ by Lemma \ref{lem:first_betti_num}.  Also
by the Novikov Signature Theorem,
$$\sigma(E(n,d,k,g))=\sigma(E(n)) + \sigma(B(d,k,g;0)) = -8n.$$
We need to check that $\kappa(E(n,d,k,g))=1$.  Since the
Poincar\'{e} dual of the canonical bundle $PD(K_{E(n)}) = (n-2)T'$
and the Poincar\'{e} dual of the canonical bundle
$PD(K_{B(d,k,g,0)})=(2g-2)T$, the fiber sum $E(n,d,k,g)$ has
canonical class $PD(K) = (n-2+2g)T$.  Therefore $K^2=0$ and
$K\cdot [\tilde{\omega}] = n-2 +2g
>0$ when $n>1$ and $g>0$.

To finish the proof of Theorem \ref{thm:main_thm} we need a
symplectic 4--manifold with $\kappa=1$ for the admissible triple
$(a,b,c)$, where $a<0$. Since $0\leq c\leq b$ with $b-c$ even, set
$k=(b+c)/2$ and note that $0\leq c\leq k$ by the fact that $0\leq
2c\leq b+c$. Then $E(-a/8, c, k,g)$ is the desired manifold for
$g\geq k$, after observing the next lemma.

\begin{lemma}
The degeneracy of $E(n,d,k,g)$ is $d$.
\end{lemma}

\begin{proof} It is straight forward to compute the ring structure
of $B(d,k,g,0)$ given that it is a product of a three manifold
with $S^1$ (See Lemma \ref{lem:degeneracy} for an example of such calculations).
\end{proof}

For symplectic manifolds with signature $-8$, we fibersum with
Dolgachev surfaces $E(1)_{p,q}$ instead of $E(1)$ to get a minimal
symplectic 4--manifold.

\section{Null admissible triples}

We do not know how much the degeneracy of a symplectic manifold $(M,\omega)$ depends on the symplectic form $\omega$.
It may be true that $d(M,\omega) \not= d(M,\omega')$ for two different symplectic forms $\omega$ and $\omega'$ on the same manifold.  One could look at another invariant, called the {\em nullity of $M$}, which only depends on the ring structure of $M$ and not on the particular symplectic form chosen.

\begin{defn} For any $\alpha\in H^1(M;{\bf R})$, and $i=1, 2$, consider
the map
$$i_{\alpha}=\cup \alpha:H^i(M;{\bf R})\longrightarrow H^{i+1}(M;{\bf R}).$$
The dimension of the linear space $\{\alpha|i_{\alpha}=0\}$ is
called the {\em $i$--nullity of $M$}, denoted $n_i(M)$.
\end{defn}

\begin{lemma} $n_1(M)=n_2(M)$. \label{lem:nullity}
\end{lemma}

\begin{proof} If $1_{\alpha}=0$, then we claim that $2_\alpha=0$.
Otherwise there is a $\gamma \in H^2(M;{\bf R})$ such that
$\alpha\cup \gamma$ is nonzero in $H^3(M;{\bf R})$. By the
Poincar\'e Duality, there is a class $\beta\in H^1(M;{\bf R})$
such that $(\alpha\cup \gamma)\cup \beta\ne 0$. This implies that
$1_{\alpha}(\beta)=\alpha\cup \beta\ne 0$, which is a
contradiction. Similarly we can prove that $2_{\alpha}=0$ implies
that $1_{\alpha}=0$.
\end{proof}

Thus we can speak simply of nullity of $M$, $n(M)$.

It follows from Lemma \ref{lem:nullity}  that the nullity is a lower bound for the degeneracy of $M$, i.e.,  $d(M,\omega)\geq
n(M)$.  From  Hodge theory we get that \Kahler\ surfaces are
of Lefschetz type, and hence have nullity zero.

One can also talk about triples $(a,b,c)$ where $c$ is the nullity of a Kodaira dimension 1 symplectic manifold (or for a symplectic manifold in general).

\begin{defn} Any triple $(a,b,c) \in \BZ^3$ is called {\em null admissible}
if  $a=8k$ where $k$ is a non-positive number, $0\leq c \leq b$, $c\not=b-1$,
and $b \geq \max\{0, 2 + a/4\}$.
\end{defn}

Observe that we require that $c\not=b-1$, for if the nullity of a manifold $M$ was $b_1(M)-1$, there would be an element in $H^1(M;\BZ)$ whose cup product square would not be zero.  Note that we no longer have a reason to require that $b-c$ be even.  In fact, the next lemma shows that
$b-c$ may be even or odd.

\begin{lem} The nullity $n(B(d,k,g;e))$ is
$$n(B(d,k,g;e)) = \left\{\begin{array}{c@{ \ \ \ }l} 0 & e=0\\
d & e\not=0, d\not=k\\
d+1 & e\not=0, d=k\end{array} \right.$$\label{lem:nullityofB}
\end{lem}

\begin{proof}
Compute the ring structure for cup products on $H^1(B(d,k,g;e); \BR)$ using the explicit basis given in Equation \eqref{eq:explicit_basis} and in the proof of Lemma \ref{lem:degeneracy}.
\end{proof}

We can use Lemma \ref{lem:nullityofB} to find null admissible triples for Kodaira dimension 1.  For example, when $a=0$ and $b=2$ (recall that $b$ can not equal 1), the lemma shows that $B(0,0,g;0)$ has nullity zero and that $B(1,1,g;1)$  has nullity 2; the triple $(0,2,1)$ is not a null admissible triple.

For the case when $a=0$ and $b=3$, we can recognize the null admissible triples $(0,3,0)$ and $(0,3,3)$ using bundle manifolds $B(1,1,g;0)$ and $B(2,2,g;1)$ respectively.  Lemma \ref{lem:nullityofB} does not give an example for the null admissible triple $(0,3,1)$. The ring multiplication for such a manifold would have to look like
\[
\begin{array}{c|ccc} & a& b& c\\ \hline
a & 0 & a\cup b & 0\\
b& b\cup a & 0 & 0\\
c & 0 & 0 & 0\end{array}
\]
for a basis $\langle a,b,c\rangle$ of $H^1(M;\BZ)$ where $a\cup b\not=0$.  We end this report with the following question.

\begin{quest} Does there exist a symplectic 4--manifold $(M,\omega)$ with $\sigma(M)=0$, $b_1(M)=3$, and nullity $n(M)=1$?
\end{quest}

{\bf Acknowledgements}\qua
The first author is partially supported by NSF grant DMS-0406021. The second author is partially supported by NSF and the
McKnight fellowship.

\Addresses\recd

\end{document}

%% file: agtout.tex
\def\ifplaintex{\expandafter\ifx\csname documentclass\endcsname\relax}

\def\gtp{{\mathsurround=0pt\it $\cal G\mskip-2mu$eometry \&\ 
$\cal T\!\!$opology $\cal P\!$ublications}}  % GT publications

\def\recd{{\small Received:\qua\receiveddate\ifx\reviseddate\relax
\else\qquad Revised:\qua\reviseddate\fi\par}} 

\def\lognumber#1{\def\thelognumber{#1}}
\def\volumenumber#1{\def\thevolumenumber{#1}}
\def\volumeyear#1{\def\thevolumeyear{#1}}
\def\papernumber#1{\def\thepapernumber{#1}}
\def\pagenumbers#1#2{\def\startpage{#1}\def\finishpage{#2}}
\def\published#1{\def\publishdate{#1}}

\def\received#1{\def\receiveddate{#1}}
\def\revised#1{\def\reviseddate{#1}}
\def\accepted#1{\def\accepteddate{#1}}
\def\asciititle#1{\def\theasciititle{#1}}

\def\asciiaddress#1{\def\theasciiaddress{#1}}
\def\asciiemail#1{\def\theasciiemail{#1}}

\long\def\asciiabstract#1{\long\def\theasciiabstract{#1}}

\let\\\par\let\thelognumber\relax\let\thevolumenumber\relax
\let\thepapernumber\relax\let\thevolumeyear\relax\let\startpage\relax
\let\finishpage\relax\let\publishdate\relax\let\receiveddate\relax
\let\reviseddate\relax\let\accepteddate\relax\let\theasciititle\relax
\let\theasciiauthors\relax\let\theasciiaddress\relax
\let\theasciiabstract\relax

\let\theasciiemail\relax

\ifplaintex
\font\logobig=cmssbx10 scaled 3836
\font\logomed=cmssbx10 scaled 2557
\else
\font\logobig=cmssbx10 scaled 4200
\font\logomed=cmssbx10 scaled 2800
\fi

\long\def\makeagttitle{   %%% start of definition of \makeagttitle
\count0=\startpage
\agt\hfill      %   Journal title (top left) 
\hbox to 45truept{\vbox to 0pt{\vglue -13truept{\logomed A\kern -.37em{\logobig 
T}\kern -.38em G}\vss}\hss}
\break
{\small Volume \thevolumenumber\ (\thevolumeyear)
\startpage--\finishpage\nl
Published: \publishdate}

\vglue .25truein

{\parskip=0pt\leftskip 0pt plus
1fil\def\\{\par\smallskip}{\Large\bf\thetitle}\par\medskip} \vglue
0.05truein

{\parskip=0pt\leftskip 0pt plus 1fil\def\\{\par}{\sc\theauthors}
\par\medskip}%
 
\vglue 0.03truein 

{\small\leftskip 25truept\rightskip 25truept{\bf Abstract}\stdspace\theabstract

{\bf AMS Classification}\stdspace\theprimaryclass
\ifx\thesecondaryclass\relax\else; \thesecondaryclass\fi\par
{\bf Keywords}\stdspace \thekeywords\par}\vglue 7truept

}   %%%% end of definition of \makeagttitle

\ifplaintex
\font\phead=cmsl9 scaled 950
\font\pnum=cmbx10 scaled 913
\font\pfoot=cmsl9 scaled 950
\headline{\vbox to 0pt{\vskip -4.5mm\line{\small\phead\ifnum
\count0=\startpage ISSN 1472-2739 (on-line) 1472-2747 (printed)
\hfill {\pnum\folio}\else\ifodd\count0\def\\{ }% 
\ifx\theshorttitle\relax\thetitle\else\theshorttitle\fi\hfill{\pnum\folio}
\else\def\\{ and }{\pnum\folio}\hfill\ifx\theshortauthors\relax\theauthors
\else\theshortauthors\fi\fi\fi}\vss}}
\footline{\vbox to 0pt{\vglue 0mm\line{\small\pfoot\ifnum\count0=\startpage
\copyright\ \gtp\hfill\else
\agt, Volume \thevolumenumber\ (\thevolumeyear)\hfill\fi}\vss}}
\else
\font=cmsl9 scaled 1050
\font\lnum=cmbx10 
\font=cmsl9 scaled 1050
\makeatletter
\def\@oddhead{{\small\lhead\ifnum\count0=\startpage ISSN 1472-2739 
(on-line) 1472-2747 (printed)\hfill {\lnum\number\count0}\else\ifodd\count0
\def\\{ }\ifx\theshorttitle\relax \thetitle \else\theshorttitle\fi\hfill
{\lnum\number\count0}\else\def\\{ and }{\lnum\number\count0}
\hfill\ifx\theshortauthors\relax 
\theauthors\else\theshortauthors\fi\fi\fi}}\def\@evenhead{\@oddhead}
\def\@oddfoot{\small\lfoot\ifnum\count0=\startpage\copyright\ \gtp\hfill\else
\agt, Volume \thevolumenumber\ (\thevolumeyear)\hfill\fi}
\def\@evenfoot{\@oddfoot}
\makeatother
\fi
\let\maketitlepage\makeagttitle
\let\makeshorttitle\maketitlepage
\let\maketitle\maketitlepage

\newwrite\gtoutfile
\long\gdef\makeheadfile{  %%% start of definition of \makeheadfile
{\def\\{, }\def\s{ }
\immediate\openout\gtoutfile head.xxx
\immediate\write\gtoutfile{Proxy-for: \ifx\theasciiauthors\relax
\theauthors\else\theasciiauthors\fi\s<\ifx\theasciiemail\relax\theemail\else\theasciiemail\fi>}
\immediate\write\gtoutfile{\noexpand\\}
\immediate\write\gtoutfile{Authors: \ifx\theasciiauthors\relax
\theauthors\else\theasciiauthors\fi}
{\def\\{ }\immediate\write\gtoutfile{Title: \ifx\theasciititle\relax
\thetitle\else\theasciititle\fi}}
\immediate\write\gtoutfile{Subj-class: GT or SG, GR etc}
\immediate\write\gtoutfile{MSC-class: \theprimaryclass\ifx\thesecondaryclass\relax\else, \thesecondaryclass\fi}
\immediate\write\gtoutfile{Journal-ref: Algebr. Geom. Topol. \thevolumenumber\s
(\thevolumeyear) \startpage-\finishpage}
\immediate\write\gtoutfile{Comments: Published by Algebraic and
Geometric Topology at}
\immediate\write\gtoutfile{\s\s\s  http://www.maths.warwick.ac.uk/agt/AGTVol\thevolumenumber/agt-\thevolumenumber-\thepapernumber.abs.html}
\immediate\write\gtoutfile{\noexpand\\}
\immediate\write\gtoutfile{}
\ifx\theasciiabstract\relax
\immediate\write\gtoutfile{\theabstract}\else
\immediate\write\gtoutfile{\theasciiabstract}\fi
\immediate\write\gtoutfile{}
\immediate\write\gtoutfile{\noexpand\\}
\immediate\write\gtoutfile{}
\immediate\closeout\gtoutfile}}  %%% end of definition of \makeheadfile

\def\maketitlepage{\makeagttitle\makeheadfile}
\let\makeshorttitle\maketitlepage
\let\maketitle\maketitlepage